\documentclass[12pt,a4paper]{amsart}

\usepackage{amsmath,amsthm,amssymb}
\usepackage{amscd}
\newtheorem{thm}{Theorem}[section]

\newtheorem{prop}[thm]{Proposition}
\newtheorem{cor}[thm]{Corollary}

\newtheorem{theorem}{Theorem}
\newtheorem{lemma}{Lemma}[section]

\newcommand{\Res}{\textrm{Res}}
\newcommand{\beq}[1]{\begin{equation}\label{#1}}
\newcommand{\eeq}{\end{equation}}
\newcommand{\bep}{\begin{proof}}
\newcommand{\eep}{\end{proof}}
\newcommand{\bec}[1]{\begin{cor}\label{#1}}
\newcommand{\eec}{\end{cor}}
\newcommand{\bepr}[1]{\begin{prop}\label{#1}}
\newcommand{\eepr}{\end{prop}}
\newcommand{\bel}[1]{\begin{lemma}\label{#1}}
\newcommand{\eel}{\end{lemma}}
\newcommand{\bet}[1]{\begin{theorem}\label{#1}}
\newcommand{\eet}{\end{theorem}}

\newcommand{\eps}{\ensuremath{\varepsilon}}

\newcommand{\bN}{\ensuremath{\mathbb{N}}}

\newcommand{\h}{{\rm h}}

\begin{document}


 \title{Searching for Diophantine quintuples}

\author{Mihai Cipu and Tim Trudgian}
\date{\today}

\begin{abstract}
We consider Diophantine quintuples $\{a, b, c, d, e\}$. These are 
sets of distinct positive integers, the product of any two elements 
of which is one less than a perfect square. It is conjectured that 
there are no Diophantine quintuples; we improve on current estimates 
to show that there are at most $1.18\cdot 10^{27}$ Diophantine quintuples.
 
\end{abstract}

\maketitle 

\section{Introduction} \label{sec1}
Define a Diophantine $m$-tuple as a set of $m$ integers $\{a_{1}, \ldots, a_{m}\}$ 
with $a_{1} < a_{2} < \ldots < a_{m}$,  such that $a_{i} a_{j} +1$ 
is a perfect square for all $1\leq i<j\leq m$. Throughout the rest 
of this article we frequently refer to $m$-tuples, and not to Diophantine 
$m$-tuples.

It is conjectured that there are no quintuples --- see~\cite{ahs,gib}. 
Successive authors (see, e.g., Table 1 in \cite{tim}) have reduced 
the bound on the possible number of quintuples. The best such 
published bound is $2.3\cdot 10^{29}$ by Trudgian \cite{tim}. The 
purpose of this paper is to improve on this in the following theorem.

\begin{theorem}\label{prince}
There are at most $1.18\cdot 10^{27}$ Diophantine quintuples.
\end{theorem}

In \S \ref{sec2} we collect some ancillary results that aid the computational 
search for quintuples. In \S \ref{sec3} we obtain bounds on the relative 
sizes of elements in a quintuple. We use this in \S \ref{sec4} with results 
on linear forms of logarithms to obtain upper bounds on the second-largest 
element in a quintuple. In \S \ref{sec5} we examine some number-theoretic sums, 
which enable us to bound the total number of quintuples. We present two new
arguments in \S \ref{sec6} that enable us to make a further saving, and 
ultimately to prove Theorem \ref{prince}.

We are grateful to Adrian Dudek who suggested the asymptotic form in (\ref{hunger}).

\section{Discards} \label{sec2}

It is known that every triple $\{a, b, c\}$ can be extended to a 
quadruple of a certain form. This is dubbed the `regular' 
quadruple and is denoted as $\{a, b, c, d_{+}\}$. If a double or 
a triple cannot be extended to a non-regular quadruple, then it 
cannot be extended to a quintuple. We call such doubles or triples 
\textit{discards}. The doubles $\{k, k+2\}$ \cite{nou} (see also 
\cite{bdm})  are discards for $k\geq 1$. For an extensive list of 
discards, one may see \cite[\S 2.1]{tim}. The following result 
allows us to recognise many discards.

\bel{lem:lex}
Let $\{a,b,c,d\}$ be a Diophantine  quadruple with
$a<b<c<d_+<d$.

 $\bullet$ If $b<2\, a$ then $b > 21000$.

 $\bullet$ If $2\, a\le b \le 12\,  a$ then  $b > 130000$.

 $\bullet$ If $b>12\,  a$ then $b>4001$.
\eel
\bep 
The only difference between this lemma and Lemma~3.4 in~\cite{cff} is the 
exclusion of the value $b=4001$ in the last case. Indeed, 
a pair $\{a,4001\}$  with $12a<4001$ cannot be extended 
because the equation $4001a+1=r^2$ has unique integer 
solution $r<4001$, namely $r=4000$, which entails $a=3999$. 
\eep 

\bel{lem:t11} \emph{(\cite[Theorem~1.1]{cff})}
Let  $\{a,b,c,d,e\}$ be a  quintuple  with $a<b<c<d<e$  and put 
$g=\gcd(a,b)$. Then $b > 3ag$. If moreover $c>a+b+2\sqrt{ab+1}$ 
then $b > \max\{24 \, ag,2\, a^{3/2}g^2\}$.
\eel

\bel{lem:t13} \emph{(\cite[Theorem~1.3]{cff})}
Let  $\{a,b,c,d,e\}$ be a  quintuple  with $a<b<c<d<e$  and  
$c=a+b+2\sqrt{ab+1}$. Then $b < a^3$ and $\gcd(b,c)=1$.
In particular, at least one of $a$, $b$ is odd.
\eel

Examination of the relative size of  entries in a
quintuple  has the following outcome. 

\bel{lem:type}
Any quintuple $\{a,b,c,d,e\}$ with $a<b<c<d<e$ must be of one of the 
types listed below:
\begin{enumerate}
 \item[(A)]  $4a < b$ and  $ 4ab+b+a< c < b^{3/2}$,
 \item[(B)]   $4a < b$ and  $c=a+b+2\sqrt{ab+1}$,
 \item[(C)]  $4a < b$ and  $ c > b^{3/2}$,
 \item[(D)]  $b < 4a$ and  $c=a+b+2\sqrt{ab+1}$.
\end{enumerate}
\eel \bep 
The claim follows from  the previously cited results having 
in view \cite[Lemma~4.2]{eff} and Subsection 2.2 in~\cite{eu}.
\eep 

\section{Exploiting the connection with Pellian equations} \label{sec3}

The entries in a quadruple are severely restricted in that they 
appear as coefficients of three generalized Pell equations
that must have at least one common solution in positive
integers. Each component of such a solution is obtained as a 
common term of two second-order linearly recurrent sequences,
giving rise to relations of the type $z=v_m=w_n$ for some
positive integers $m$  and $n$. A key ingredient in the study 
of Diophantine sets is a  relationship between the parameters 
$m$, $n$, and the values in the set in question.

Our next result is of this kind. It improves on several versions
already in the literature --- see, e.g., \cite{eu,tim,wh}.

\bepr{pr:mlow}
Let $\{A,B,C,D\}$ be a quadruple with $A<B<C<D$ for which $v_{2m}=w_{2n}$ 
has a solution with $2n\ge m\ge n \ge 2$, $m\ge 3$. Suppose that 
$A\ge A_0$, $B\ge B_0$, $C\ge C_0$, $B> \rho A$ for some positive 
integers $A_0$, $B_0$, $C_0$, and a real number $\rho  > 1$. Then 
\[
 m > \alpha B^{-1/2}C^{1/2},
\]
where $\alpha$ is any real number satisfying both inequalities
\beq{eq:al1}
\alpha^2 +\bigl(1+\tfrac{1}{2}B_0^{-1}C_0^{-1} \bigr) \alpha\le 4,
\eeq 
\beq{eq:al2}
3\alpha^2 +\left(4B_0 (\lambda +\rho^{-1/2})+ %
2(\lambda +\rho^{1/2})C_0^{-1} \right) \alpha\le 4B_0,
\eeq
with $\lambda = (A_0+1)^{1/2} (\rho A_0+1)^{-1/2}$.

Moreover, if $C^{\tau}\ge \beta B$ for some positive real numbers 
 $\beta$  and $\tau$  then
\[
 m > \alpha \beta ^{1/2} C^{(1-\tau)/2}.
\]
\eepr \bep 
We assume that $ m \le \alpha B^{-1/2}C^{1/2}$ and aim at establishing 
a contradiction if $\alpha$ is too small.
We start from the congruence (see, e.g., \cite[Lemma~4]{dub})
\beq{eq:con}
\eps Am^2 +Sm \equiv \eps Bn^2+Tn \pmod{4C} \quad 
\mbox{for some} \quad \eps =\pm 1,
\eeq 
where $S = \sqrt{AC + 1}$ and $T=\sqrt{BC +1}$.
Since
\[
\vert Am^2-Bn^2\vert < \max \{Am^2,Bn^2 \}\le  Bm^2 \le \alpha ^2\, C
\]
and
\begin{align*}
\vert Sm-Tn \vert & < \max \{Sm,Tn \} \le Tm \le \alpha B^{-1/2}C^{1/2} 
\sqrt{BC+1} \\
& <  \alpha B^{-1/2}C^{1/2}  \bigl( B^{1/2}C^{1/2}+\tfrac{1}{2} B^{-1/2}C^{-1/2} \bigr) \le 
\alpha \bigl(1+\tfrac{1}{2}  B_0^{-1}C_0^{-1} \bigr) C,
\end{align*}
then, if $\alpha$ satisfies~\eqref{eq:al1},  the congruence~\eqref{eq:con}
becomes the equality $Am^2-Bn^2=\eps (Tn-Sm)$. Multiplication by
$Tn+Sm$ followed by rearrangements results in the equality 
\beq{eq:star}
 (Bn^2-Am^2)\bigl( C+\eps (Tn+Sm)\bigr)=m^2-n^2.
\eeq 
Note that $Bn^2=Am^2$ entails $m^2=n^2$, so that $A=B$: a contradiction.
Hence, for $m=n$ one necessarily has $C=Tn+Sm$, while for $m>n$ one finds 
that $Bn^2-Am^2$ divides the positive integer $m^2-n^2$, so that 
$m^2-n^2\ge \vert Am^2-Bn^2\vert$. This gives the following inequality
\[
\frac{m^2}{n^2} \ge \frac{B+1}{A+1}.
\]
Having in view the lower bounds for $A$ and $B$, we obtain
\[
 \frac{m^2}{n^2} > \frac{\rho A+1}{A+1} \ge \frac{\rho A_0+1}{A_0+1}=
  \frac{1}{\lambda^2}. 
\]
From~\eqref{eq:star}, $m\le 2n $, and the definitions of $S$ and $T$, 
we conclude that 
\begin{align*}
 C & \le Tn+Sm + m^2-n^2 < \lambda m\sqrt{BC+1}+m\sqrt{AC+1} + 
 \tfrac{3}{4}\, m^2 \\
 & \le \tfrac{3}{4}\,  \alpha^2 B^{-1}C + \alpha B^{-1/2}C^{1/2} 
 \left( \lambda \sqrt{BC+1}+\sqrt{ \rho^{-1} BC+1} \right)  \\
  & < \tfrac{3}{4}\,  \alpha^2 B^{-1}C + \alpha C
 \left( \lambda \bigl(1+\tfrac{1}{2}  B^{-1}C^{-1} \bigr) + \rho^{-1/2} %
 \bigl(1+\tfrac{1}{2}   \rho  B^{-1}C^{-1} \bigr) \right)  \\
 & \le \tfrac{3}{4}\,  \alpha^2 B_0^{-1}C + \alpha C
 \left( \lambda \bigl(1+\tfrac{1}{2}  B_0^{-1}C_0^{-1} \bigr) + \rho^{-1/2} %
 \bigl(1+\tfrac{1}{2}  \rho  B_0^{-1}C_0^{-1} \bigr) \right).
\end{align*}
The last expression is at most $C$ if $\alpha$ satisfies the 
inequality~\eqref{eq:al2}, whence the first inequality in the 
conclusion of our proposition. The second one is readily
obtained from what we have just proved and the hypothesis
$C^{\tau}\ge \beta B$.\eep 

\bel{lem:pell}
If $\{a,b,c,d,e\}$ is a   quintuple with $a<b<c<d<e$ then the 
following bounds for $m$ hold:

\emph{(A)}  $m> 3.3022 d^{1/4}$, \quad 
\emph{(B)}   $m> 1.5002  d^{2/7}$, 
\quad \emph{(C)}  $m> 2.0604 d^{3/10}$,

 \emph{(D)}  $m> 1.0080 d^{1/3}$.
\eel \bep 
This is an application of the result just proved for
$(A,B,C)=(a,b,d)$ in cases (A)--(C) and for $(A,B,C)=(a,c,d)$ 
in the remaining case.
We use Proposition~\ref{pr:mlow}  with carefully chosen values
for parameters in ranges suggested by 
Lemmas~\ref{lem:lex}--\ref{lem:t13}. In~\cite{fu1} it was shown 
that the first hypothesis required in Proposition~\ref{pr:mlow} 
is satisfied. It is also known that one has $d> 4abc+a+b+c$
(see, for instance, the proof of Lemma~6 in~\cite{dus}).

In case (A) Lemma~\ref{lem:t11} hints to consider separately
values of $a$ less than 144 since then one has 
$B=b>\max \{24 a, 2a^{3/2}\}=24a=24A$. A short computer
search finds that $B_0=4095$ and $\frac{b}{a}\ge 
\frac{21001}{143}>146=:\rho$. Clearly we must put $A_0=1$. From
\[
C = d> 4abc+a+b+c>(4ab+1)(4ab+a+b) > (16a^2+4a)b^2
\]
it follows that  $\tau=1/2$, $\beta =(16A_0^2+4A_0)^{1/2}$,
$C_0>3.35\cdot 10^8$ are admissible choices. Both 
inequalities~\eqref{eq:al1} and~\eqref{eq:al2} are satisfied 
by $\alpha =1.56155$.

Still in case (A), when $a\ge 144$ one puts $A_0=144$, $B_0=4002$
(by Lemma~\ref{lem:lex}), $\rho=24$ (see Lemma~\ref{lem:t11}),
 $\tau=1/2$, $\beta =(16A_0^2+4A_0)^{1/2}$, whence
$C_0>5.32\cdot 10^{12}$ and  $\alpha =1.56155$.

Having in view  Lemma~\ref{lem:lex}, in case (B) we first examine 
the subcase $4a < b \le 12 a$. Then $B_0=130001$, which implies
$A_0=10834$ and $\rho=4$. From
\[
 c > b \bigl(1+12^{-1}+2\cdot 12^{-1/2}\bigr)=
 (1+12^{-1/2})^2B
\]
and  $a^3>b$ it follows that
\[
C = d> (4ab+1)(a+b+2r) > 4(1+12^{-1/2})^2ab^2 > 4 (1+12^{-1/2})^2 B^{7/3},
\]
so that $\tau=3/7$, $\beta =(2+3^{-1/2})^{6/7}$, 
$C_0=5.68\cdot 10^{12}$. For these choices it is readily obtained that
$\alpha  = 0.9999$ is permissible.

The other possibility in case (B) is to have $b > 12a$.
Convenient values of parameters are  $\rho=12$, $A_0=16$
(from $a^3 >b > 4000$), 
$B_0=4002$, $\tau=3/7$, $\beta =2^{6/7}$,  $C_0=1.01\cdot 10^9$.
for which the same value $\alpha = 0.9999$ works.

Case (C) is similar to case (A). Now, for $a\le 143$ we see that
we can take  $A_0=1$, $B_0=4004$,  $\rho=28$. As
\[
C > 4abc > 4ab^{5/2} > 4.05\cdot 10^{9}=: C_0,
\]
we further get $\tau=2/5$,  $\beta =4^{2/5}$, whence  again  
$\alpha = 1.5615$. In the complementary subcase $a\ge 144$,
admissible values are  $A_0=144$, $B_0=4002$, $\rho=24$, 
$\tau=2/5$,  $\beta =576^{2/5}$, $C_0=5.83\cdot 10^{11}$.
Plugging these specializations into  Proposition~\ref{pr:mlow}, 
we obtain the same value for $\alpha$.

Finally, in case (D) we have $A=a <b/3$, 
$ B=c=a+b+2\sqrt{ab+1}  > (1+3^{1/2})^2A$, 
$ B \le a+b+2\sqrt{3^{-1}b(b-1)+1} < (1+3^{-1/2})^2b$,
and
\[
 C=d > 4abc > b^2c > (1+3^{-1/2})^{-4} B^3. 
\]
Therefore, $\rho=(1+3^{1/2})^2$, $\tau=1/3$, and 
$\beta=(1+3^{-1/2})^{-4/3}$.
From $130001\le b <4a$, we have $A_0= 32501$, whence $B_{0} >292504$, and
$C_0>4.04\cdot 10^{15}$.
From~\eqref{eq:al1} and~\eqref{eq:al2} we 
obtain $\alpha= 1.3660$. 
\eep 

For future reference, the values used in the 
previous proof are given in Table~\ref{tab1}. 

\begin{table}[h]
\begin{center}
\begin{tabular}{l|rrrrrr}
 Type & $A_0$ \phantom{1}  & $B_0$ \phantom{12} & $C_0$  \phantom{123}
 & $\rho$  \phantom{12} & $\beta$ \phantom{1234} &  $\tau$  \\
 \hline 
 (AI) & $1$ & 4095 & $3.35\cdot 10^8$ & 146 & $20^{1/2}$ & $1/2$\\
 (AII) & $144$ & 4002 &  $5.32\cdot 10^{12}$ & 24 & $24\cdot 577^{1/2}$ & $1/2$\\
 (BI) & 10834 & 130001 &  $5.68\cdot 10^{12}$ & 4 & $(2+3^{-1/2})^{6/7}$ &  $3/7$\\
 (BII) & 16 & 4002 & $1.01\cdot 10^9$ & 12 & $2^{6/7}$ & $3/7$\\
 (CI) & 1 & 4004 & $4.05\cdot 10^9$ & 28 & $4^{2/5}$ & $2/5$\\
 (CII) & 144 & 4002 & $5.83\cdot 10^{11}$ & 24 & $576^{2/5}$ & $2/5$ \\
 (D) & 32501 & 292504 & $4.04\cdot 10^{15}$ & $(1+3^{1/2})^2$ &   $(1+3^{-1/2})^{-4/3}$ & $1/3$ 
 \\
\end{tabular}
 \caption{Parameter values for various types of Diophantine quintuples.} \label{tab1}
\end{center}
\end{table}

The values of $\alpha$, and hence the bounds on $m$ in Lemma \ref{lem:pell}, rely on the computational bounds in Lemma \ref{lem:lex}. While it is tempting to extend these computations, such an extension would have almost no effect on the values of $\alpha$. Consider, for example, case (A): sending $B_{0}, C_{0}$ to infinity in (\ref{eq:al1}) gives $\alpha^{2} + \alpha \leq 4$. Therefore the optimal value of $\alpha$ is $1.5615528\ldots$, whereas we have $\alpha = 1.56155$. Likewise, in case (D) the optimal value is $\frac{1}{2}(1 + \sqrt{3}) = 1.366025\ldots$, whereas we have $1.3660$. It seems that a new idea is needed to improve substantially on the lower bounds on $m$.

\section{Employing linear forms in the logarithm} \label{sec4}

The lower bounds for the index $m$ given in the previous section 
can be complemented by inequalities derived from upper 
bounds for linear forms of logarithms of algebraic numbers.
To this end, we apply the best result of which we are aware.
\begin{thm}[Aleksentsev]\label{Alex}
Let $\Lambda$ be a linear form in logarithms of $n$ multiplicatively 
independent totally real algebraic numbers $\alpha_{1}, \ldots \alpha_{n}$, 
with rational coefficients $b_{1}, \ldots, b_{n}$. Let $h(\alpha_{j})$ 
denote the absolute logarithmic height of $\alpha_{j}$ for $1\leq j \leq n$. 
Let $d$ be the degree of the number field 
$\mathcal{K} = \mathcal{Q}(\alpha_{1}, \ldots, \alpha_{n})$, and let 
$A_{j} = \max(d h(\alpha_{j}), |\log \alpha_{j}|, 1)$. Finally, let
\begin{equation}\label{tail}
E = \max\left( \max_{1 \leq i, j \leq n} \left\{ \frac{|b_{i}|}{A_{j}} + \frac{|b_{j}|}{A_{i}}\right\}, 3\right).
\end{equation}
Then
\begin{equation*}\label{kidney}
\log |\Lambda| \geq - 5.3  n^{-n+1/2} (n+1)^{n+1}(n+8)^{2}(n+5)(31.44)^{n} 
d^{2} (\log E) A_{1}\cdots A_{n} \log(3nd).
\end{equation*}
\end{thm}
We have used the first displayed equation on \cite[p.\ 2]{alex} 
to define $E$ in (\ref{tail}): this makes our application easier.
We apply Theorem \ref{Alex} for $d=4, n=3$ and to 
\[
 \Lambda = j\log\alpha_{1} - k\log\alpha_{2} + \log\alpha_{3},
\]
with
\begin{equation*}\label{cheek}
\alpha_{1} = S+ \sqrt{AC}, \quad \alpha_{2} = T + \sqrt{BC}, \quad 
\alpha_{3} = \frac{\sqrt{B}(\sqrt{C} \pm \sqrt{A})}{\sqrt{A} \left( \sqrt{C} \pm \sqrt{B}\right)}, 
\end{equation*}
where the signs coincide.

For our purposes we do not need the exact values of
$A_j$ and $E$ as defined in Theorem \ref{Alex}:
decent estimates will suffice. 
To find these estimates we proceed as follows, keeping the notation 
and hypotheses of Proposition~\ref{pr:mlow} and supposing additionally
that $C\le C_1$ for a certain integer $C_1$.
 
We begin by noting that one has
\begin{align*}
2\log \alpha _1 & < \log (4AC+4)\le \log \bigl(4\rho ^{-1} (B-1)C+4 \bigr)
< \log (4\rho ^{-1} BC) \\
& < \log (4\rho ^{-1} \beta ^{-1}C^{1+\tau})
\end{align*}
provided that $\rho A \le B-1$. This clearly follows from $\rho A < B$
when $\rho$ is integer, as in cases (A)--(C). In case (D) we have $b\ge 3a+1$, so that 
(cf.\ the proof of Lemma~\ref{lem:pell})
\[
B=c=a+b+2\sqrt{ab+1} > 1+(1+ 3^{1/2})^{2}a=1+\rho A. 
\]

In each of the cases (A)--(D) we have $\beta \rho >4$, whence 
\begin{equation*}\label{eq:up1}
A_1 < g_1(\beta, \rho, \tau, C_1) \log C,
\end{equation*}
with
\[
g_1(\beta, \rho, \tau, C_1) :=1+\tau +
\frac{\log 4- \log (\beta \rho)}{\log C_1}. 
\]

We readily obtain the following lower bound on $A_{1}$
\begin{equation*}\label{eq:do1}
A_1 > \log (4AC)>  g_2(A_0, C_1) \log C,
\end{equation*}
with
\[
g_2(A_0, C_1) :=1+\frac{\log 4+ \log A_0}{\log C_1}. 
\]

Similar relations hold for $A_2$, namely
\[
2\log \alpha _2  < \log (4BC+4)  < \log (4 \beta ^{-1}C^{1+\tau}+4),
\]
which implies the upper bound
\[
 A_2 < g_3(\beta,  \tau, e) \log C,
\]
where
\[
 g_3(\beta,  \tau, e) :=1+\tau +
 \frac{\log 4+ \log (\beta^{-1}+ e^{-1-\tau})}{\log e},
\]
and $e=C_0$ in the cases (B), (CI), and (D) (when $\beta < 4$)
and $e=C_1$ in the remaining cases (A) and (CII).
An easily-derived lower bound for $A_2$ is 
\[
 A_2> g_4(B_0, C_1) \log C,
\]
with 
\[
 g_4(B_0, C_1) := 1+\frac{\log 4+\log B_0}{\log C_1}. 
\]

The inequalities
\[
 \frac{\sqrt{B}}{\sqrt{A}} \cdot 
  \frac{\sqrt{C}+ \sqrt{A}}{\sqrt{C}- \sqrt{B}} >
  \frac{\sqrt{B}}{\sqrt{A}} \cdot 
  \frac{\sqrt{C} +\sqrt{A}}{\sqrt{C}+ \sqrt{B}} > 1, \quad 
 \frac{\sqrt{B}}{\sqrt{A}} \cdot 
  \frac{\sqrt{C}- \sqrt{A}}{\sqrt{C}- \sqrt{B}} >1    
\]
are obvious. The modulus of the fourth algebraic conjugate of 
$\alpha _3$ is also greater than $1$ precisely when 
$\sqrt{C}(\sqrt{B}- \sqrt{A}) > 2\sqrt{AB}$. This inequality 
holds whenever
\beq{eq:con4}
 \rho B_0^{1-\tau} (\rho^{1/2}-1)^{2\tau}> 2^{2\tau}.
\eeq 
It is easy to check that (\ref{eq:con4}) is satisfied 
in each of the cases (A)--(D). One now obtains
\[
A_3= 4\h (\alpha _3) = \log \left(\frac{B^2 (C-A)^2}{g} \right),
\]
where $g$ is the content of the polynomial 
$A^2(C-B)^2X^4+4A^2B(C-B)X^3+2AB(3AB-AC-BC-C^2)X^2 
+4AB^2(C-A)X + B^2(C-A)^2$. Since $g$ is  at most the smallest
of the coefficients, which is $4A^2B(C-B)$, one has
\[
 \log \left(\frac{B (C-A)^2}{4A^2 (C-B)} \right) \le A_3 \le
 \log \bigl(B^2 (C-A)^2 \bigr).
\]

Note that $ B(C-A)<\beta^{-1}C^{1+\tau}$ readily implies  
\[
A_3 < g_5(\beta, \tau, f) \log C,
\]
with
\[
g_5(\beta, \tau, f) :=2+2\tau -\frac{2\log \beta}{\log f} 
\]
and $f=C_1$ if $\beta > 1$ and $f=C_0$ if $\beta < 1$.
A lower bound for $A_3$ is obtained with the help of the inequalities
\begin{equation*}\label{eq:do3}
A_3  \ge \log \left(\frac{B (C-A)^2}{4A^2 (C-B)} \right) > \log  
\left(\frac{\beta \rho^2 C^{1-\tau} (1-A_0 C_1^{-1})^2}{4 (1-\rho A_0 C_1^{-1})} \right),
\end{equation*}
which entail 
\[
A_3 > g_6(\beta,\rho,\tau, A_0, C_1) \log C,
\]
where
\[
g_6(\beta,\rho,\tau, A_0,C_1):= 1-\tau +
\frac{\log (\tfrac{1}{4}\beta\rho^2)+2\log (1-A_0C_1^{-1})-\log (1-4C_1^{-1})}{\log C_1}.
\]

On noting that for all relevant values of parameters one has
$g_2(A_0,C_1) < g_4(B_0,C_1)$ and using the  inequality 
$g_2 > g_6$ (which follows, for $C_1>10^{12}$, from 
$16C_1^{2\tau}(1-\beta^{-1}) > \beta\rho^2$ if  $\beta > 1$ and  from
$16 (1+3^{-1/2})^{8/3} C_1^{2/3} (1-\rho A_0C_1^{-1}) > \rho ^3$
in  case (D))
as well as the known relation $k\le j$, we find that we may take 
\[
 E\le \frac{2j}{g_6(\beta,\rho,\tau,A_0,C_1) \log C_0}.
 \]
Hence, Theorem \ref{Alex} yields the following corollary.
\bec{eq:aleup}
\[
-\log \Lambda \le 1.5013\cdot 10^{11}g_3g_5 (2\log \alpha _1) 
(\log ^2 C) \log \left( \frac{2j}{g_6 \log C_0}\right).
\]
\eec 
Corollary~\ref{eq:aleup} bounds $\Lambda$ from below; we can bound 
$\Lambda$ from above using
 Eq.~(4.1) in~\cite{fu1}, which states that
\[
 0<  \Lambda < \frac{8}{3} AC \alpha _1^{-2j}.
\]
Comparison with Corollary~\ref{eq:aleup} gives the main result 
of this section. 

\bepr{pr:3log}
\[
j < 1.50131\cdot 10^{11}g_3g_5 (\log ^2 C) 
\log \left( \frac{2j}{g_6 \log C_0}\right).
\]
\eepr 

Set  $j=2m$ in Proposition~\ref{pr:3log} and use Lemma~\ref{lem:pell}  
with the values given in Table~\ref{tab1}
and $C_1=10^{72.188}$ in all cases,
as per~\cite[Theorem~1.2]{eu}. 
We thus get 
a new upper bound on $d$ that we take as $C_1$ in a new iteration 
of this procedure. Slightly better bounds result by taking much 
higher $C_0$ (just below the value for $C_1$ considered in the 
same iteration). We record our computations in the following theorem.

\bet{th:dup}   
If $\{a,b,c,d,e\}$ is a   quintuple with $a<b<c<d<e$ then the 
following bounds for $d$ hold:                  

\emph{(A)}  $d< 10^{67.859} < 7.228 \cdot 10^{67}$, \quad 
\emph{(B)}   $d< 10^{60.057} < 1.141  \cdot 10^{60}$,         

 \emph{(C)}   $d< 10^{56.528} <  3.373 \cdot 10^{56}$, \quad
 \emph{(D)}   $d< 10^{51.416} <  2.603  \cdot 10^{51}$.
\eet 
 
We close this section with a remark concerning the size of
the smallest entry in a quintuple. Although it has no immediate 
bearing on the next section, further improvements on $d$ should 
enable future researchers to enumerate all possible triples. 
Recording the maximal size of $a$ should aid this goal.

\bepr{pr:acasea} 
The only quintuples that could arise from case \emph{(A)} are those
in which $a < 7.4\cdot10^{7}$.
\eepr \bep 
The triples in case (A) must satisfy $b^{3/2} > c > 4ab + b + a$, so
that, in particular $a< b^{1/2}/4$. Some quick computations give that for 
$A_{0} = 7.4\cdot10^{7}$  one obtains $d< 6.1\cdot10^{50}$. From
$d>  4abc > 16a^2b^2> (16a^2)^3$  it then follows $a< 7.29\cdot 10^{7}$,
a contradiction.
\eep

 \section{Bounding the total number of quintuples} \label{sec5}
 
In this section we combine the methods of \cite{eu} and \cite{tim} 
in bounding certain arithmetical sums. Let 
\[ 
E(x) = \sum_{n=1}^{x} 2^{\omega(n)}, \quad 
F(x) = \sum_{n=1}^{x} \frac{2^{\omega(n)}}{n}, \quad 
G(x) = \sum_{n=1}^{x} \frac{2^{\omega(2n-1)}}{2n-1},
\]
where all the sums are defined for real values of $x\geq 1$. Bounds 
for these sums can be used as in the proof of Lemma 4.4 in \cite{eu} to prove that
\begin{equation}\label{milk}
\sum_{n=2}^{N} d_{H}(n^{2} -1) \leq 2N G\left( \frac{H+1}{2}\right) + 
N G\left( \frac{H+4}{8}\right)  + N G\left( \frac{H+2}{4}\right) 
 + N G\left( \frac{H}{8}\right),
 \end{equation}
 where $d_{H}(n)$ counts the number of divisors of $n$ that do not exceed $H$. 
 The function $d_{H}(n^{2} -1)$ arises naturally when considering the number 
 of doubles $\{a, b\}$ satisfying certain restrictions.

The following lemma gives good bounds on $E(x)$ and $F(x)$.
\begin{lemma}[Lemma 13 in \cite{tim}]\label{Lem:tim}
 For all $x\geq 1$ we have
 \begin{equation}\label{actual1}
  \begin{split}
  \sum_{n\leq x} \frac{2^{\omega(n)}}{n} &\leq 3 \pi^{-2} \log^{2} x + 1.3948\log x + 0.4107+ 3.253x^{-1/3},\\
  \sum_{n\leq x} 2^{\omega(n)} &\leq 6 \pi^{-2} x \log x + 0.787x + 8.14 x^{2/3}  - 0.3762 .\\
  \end{split}
  \end{equation}
\end{lemma}
One can show, using Perron's formula and calculating residues, that
\begin{equation*}\label{help}
\begin{split}
E(x) &\sim \frac{6}{\pi^{2}} x \log x +\frac{6}{\pi^{4}}\left( \pi^{2}(2\gamma -1) - 12 \zeta'(2)\right) x,\\
F(x) &\sim \frac{3}{\pi^{2}} x \log x +\frac{12}{\pi^{4}}\left( \pi^{2}\gamma - 6 \zeta'(2)\right) x,
\end{split}
\end{equation*}
where 
\begin{equation*}\label{help2}
\frac{6}{\pi^{4}}\left( \pi^{2}(2\gamma -1) - 12 \zeta'(2)\right) = 0.78687\ldots, \quad \frac{6}{\pi^{4}}\left( \pi^{2}(2\gamma -1) - 12 \zeta'(2)\right) = 1.39479\ldots. 
\end{equation*} 
This shows that up to three decimal places, the bounds in Lemma \ref{Lem:tim} agree with the asymptotic expansions to the first two terms. 

Similarly, one can calculate the asymptotic order of $G(x)$. For $\Re(s)>1$ we have
\begin{equation*}\label{tree}
\sum_{n=1}^{\infty} \frac{2^{\omega(n)}}{n^{s}} = \frac{\zeta^{2}(s)}{\zeta(2s)} = \prod_{p} \left( 1 + \frac{2}{p^{s}} + \frac{2}{p^{2s}} + \cdots\right).
\end{equation*}
We can remove the contribution of $p=2$ to yield
\begin{equation*}\label{yam}
\sum_{n=1, n\; \textrm{odd}}^{\infty} \frac{2^{\omega(n)}}{n^{s}} =  \frac{1- 2^{-s}}{1 + 2^{-s}}\frac{\zeta^{2}(s)}{\zeta(2s)},
\end{equation*}
whence, by Perron's formula we have
\begin{equation}\label{hunger}
\begin{split}
G(x) = \sum_{n\leq 2x, n\; \textrm{odd}} \frac{2^{\omega(n)}}{n} &\sim \Res \left\{ \frac{1- 2^{-(s+1)}}{1 + 2^{-(s+1)}}\frac{\zeta^{2}(s+1)}{\zeta(2(s+1))}\frac{(2x)^{s}}{s}; \quad s=0\right\}\\
&\sim \frac{1}{\pi^{2}} \log^{2} x + \frac{2}{3 \pi^{4}} \{\pi^{2}( 6\gamma  + 7  \log 2) - 36 \zeta'(2)\} \log x.
\end{split}
\end{equation}
Were one to use this in the proof of Lemma 4.4 in \cite{eu} one would have
\begin{equation}\label{duke}
\sum_{n=2}^{N} d_{H}(n^{2} -1) \leq \frac{7N}{\pi^{2}} \log^{2} H + \ldots.
\end{equation}
It would take considerable effort to furnish an explicit version of (\ref{duke}). 
Even if the lower order terms in (\ref{duke}) were negligible, one would only
save a factor of two for the values of $H$ and $N$ used in our calculations. 
We have therefore not pursued this.\footnote{We note that Dudek \cite{DudekSum} has recently shown that $\sum_{n=2}^{N} d(n^{2} -1) \sim 6\pi^{-2} N \log^{2}N$.}
 
We content ourselves with the bound for $G(x)$ as given in Lemma 4.3 in \cite{eu}, namely
$$ G(x) < \frac{3}{2\pi^{2}} \log^{2} x + 3.1227147 \log x + 3.56851 + \frac{0.525}{x}.$$
Using this and (\ref{actual1}) in (\ref{milk}) we follow the proof of Lemma 4.4 in \cite{eu} to prove
\begin{lemma}\label{lem:final}
Let $d_{H}(n)$ denote the number of positive integers $e$ such that $e|n$ and $e\leq H$. Then, for any $N\geq 2$ and $H\geq 1$ we have
\begin{equation*}
\sum_{n=2}^{N} d_{H}(n^{2} -1) < N\left\{ \frac{9}{\pi^{2}} \log^{2} H + 11.1468 \log H - 0.957 + \frac{24 \log H}{\pi^{2}H} + \frac{44.14}{H}\right\}.
\end{equation*}
\end{lemma}

We now proceed to examine the number of quintuples that could arise 
from each of the triples (A)--(D).

\subsection{Case (A)}
This is the most damaging case in our considerations. We have 
$r< (d/16)^{1/4}$, whence, by Theorem \ref{th:dup} we have 
$r< 4.611\cdot 10^{16} = R_{A}$. Using Lemma \ref{lem:final} we find 
that the number of doubles is at most
\[
\frac{1}{2} \sum_{r=3}^{R_{A}} d_{R_{A}}(r^{2} -1) < 4.080\cdot 10^{19}.
\] 
Since $b< (d/20)^{1/2}< 1.9011 \cdot 10^{33}$ we find that $b$ could 
have as many as 23 distinct prime factors. We find that the number 
of quintuples is therefore bounded by
\begin{equation}\label{case:a}
3 \cdot 4 \cdot 2^{24} \cdot 4.08 \cdot 10^{19} \leq 8.215\cdot 10^{27}.
\end{equation}
Since the number of possible quintuples originating from case (A) 
is by far the largest, we devote \S\ref{sec6} to reducing this 
number slightly.

  \subsection{Case (B)}
Since $b>4a$ we have $b>2r$, whence $c> 4r +a$. Since $d>4abc$ this 
shows that $d> 4(r^{2} -1)(4r +2) > 16r^{3}$. From Theorem \ref{th:dup} 
we therefore have $r\leq 4.147 \cdot 10^{19} = R_B$. By Lemma 
\ref{lem:final} the number of doubles $\{a, b\}$ is at most
\[
\frac{1}{2} \sum_{r=3}^{R_B} d_{R_B}(r^{2} -1) < 4.91\cdot 10^{22}.
\]
Since there are at most four ways of extending a quadruple to a quintuple 
we find that the total number of quintuples is bounded above by
\begin{equation}\label{case:b}
2.0 \cdot 10^{23}.
\end{equation}

 \subsection{Case (C)}
We proceed as in case 2(iii) in~\cite{tim}. We consider the cases 
$a> \eta$ and $a\leq \eta$ and optimise over $\eta$. In the former case, 
we have $d> 4abc > 4\eta b^{5/2}$ so that $b< (d/(4\eta))^{2/5}:= N_{3a}$. 
Hence, by Lemma 3.3 in \cite{eff}, the number of quintuples is at most
\begin{equation}\label{quince3a}
\frac{N_{3a}}{6} \left( \log N_{3a} + 2\right)^{3}\cdot 8 \cdot 5 \cdot 4.
\end{equation}

When $a\le \eta$, we have $b< (d/(4a))^{2/5}$ so that 
$r^{2} = ab + 1 < a(d/(4a))^{2/5} +1$. Thus
\begin{equation*}
r< \sqrt{1 + \left( \frac{\eta^{3} d^{2}}{16}\right)^{1/5}}= N_{3b}.
\end{equation*}
We apply Lemma \ref{lem:final} with $H = \eta$ and $N = N_{3b}$. Since 
$b< (d/4)^{2/5} < 2.35\cdot 10^{22}$ we have $\omega(b) \leq 17$. 
Following the proof in \cite{eff} we deduce that the number of quintuples is at most
\begin{equation}\label{quince3b}
4 \cdot 2^{17} \cdot 5 \cdot 4 \cdot N_{3b}\left( \frac{9}{\pi^{2}} \log^{2} \eta + 11.1468\log \eta -0.957 + \frac{24 \log \eta}{\pi^{2} \eta} + \frac{44.14}{\eta}\right).
\end{equation}
We find that we can minimise the maximum of (\ref{quince3a}) and 
(\ref{quince3b}) at $\eta = 6.76\cdot 10^{10}$. Hence the number of 
quintuples is at most
\begin{equation}\label{case:c}
2.41\cdot 10^{22}.
\end{equation}

 \subsection{Case (D)}
We have $b< (4d/9)^{1/3}$ so that, by Theorem \ref{th:dup}, we have 
$b< 1.05\cdot 10^{17} = R_{D}$. The number of doubles $\{a, b\}$ is 
therefore bounded by $2\sum_{b=4}^{R_{D}} 2^{\omega(b)}$. We use this 
and Lemma~\ref{Lem:tim} to prove that the number of quintuples is at  
most
\begin{equation}\label{case:d}
 2.07 \cdot 10^{19}.
\end{equation}

\section{Improvements to case (A)} \label{sec6}
Here we investigate two methods. The first, in \S \ref{sec:6.1}, 
reduces the bound on $\omega(b)$ from 23 to 22, thereby saving a 
factor of 2 in the estimate recorded in~\eqref{case:a}. The second, 
in \S \ref{sec:6.2}, splits up the sum over $b$ with 
$\omega(b)$ held constant. This saves a factor of about $3.58$.

\subsection{Removing one prime factor from $b$}\label{sec:6.1}
Let $(p_n)_{n\in \bN}$ denote  the  sequence of prime numbers, 
and consider those $b$ satisfying
\begin{equation}\label{earl}
b_0:=\prod_{i=1}^{23} p_{i} \approx 2.67 \cdot 10^{32} \le b < 1.9011 
\cdot 10^{33}, \quad \omega(b) =23.
\end{equation}
We aim at enumerating all such $b$ in (\ref{earl}). We shall show that 
none of these values of $b$ can appear as the second-smallest
element of a quintuple. This then shows that  
$\omega(b) \leq 22$, and leads immediately to a saving of a 
factor of 2 in (\ref{case:a}).
 
Suppose $\{a,b,c,d,e\}$ is a quintuple. In case (A),
Theorem~\ref{th:dup} gives the bound $d< UD:=10^{67.859}$. 
When $b$ is restricted as in~\eqref{earl} we 
find that $2$ divides $b$, since, if not, the smallest  $b$ can be 
is $\prod_{i=1}^{23} p_{i} /2 \cdot p_{24} > 1.18\cdot 10^{34}$. 
Continuing in this way we find that $2, 3, 5, 7, 11$ must all divide $b$. 
 
From $  4a(4a+1)b^2 <UD$ it then follows $a\le 7$.  Moreover,
as the corresponding $r$ is odd, $ab$ is a multiple of $8$,
whence $b\equiv 0 \pmod 8$ for odd $a$ and $b\equiv 0 \pmod 4$
for $a\equiv 2 \pmod 4$. Hence, each such $b$ is obtained from
$b_1=b_1(a)$ by replacing $v$ of its factors $p_6$, $\ldots$, 
$p_{23}$ by other $v$ primes $p_{k_1}$, $\ldots$, $p_{k_v}$, 
where $24\le k_1 < \cdots < k_v$, and then multiplying by some
positive integer $q$  such that the result is at most 
\[
UB=UB(a,UD):=UD^{1/2} (16a^2+4a)^{-1/2}. 
\]
Here $b_1(a)=4b_0$ if $a$ 
is odd, $b_1(a)=2b_0$ if $a=2,6$, and $b_1(a)=b_0$ otherwise.

We now present a detailed exposition of the idea sketched above. 
All computations have been performed with GP scripts~\cite{gp}.
Clearly, the maximal $v$ is determined from the condition
\[
 \frac{ p_{24} p_{25}\cdots  p_{23+v}}{ p_{23} p_{22}\cdots  p_{24-v}}
 < \frac{UB}{b_1}.
\]
A short computer search gives $v=3$ for $a=2$ or $4$; $v=2$ 
for $a=1$; $v=0$ for the other values  $a\le 7$.

Next for each $u=1$, $2$, $\ldots$, $v$ we look for the 
largest index $K=K(u)$ satisfying
\[
 \frac{ p_{24} p_{25}\cdots  p_{22+u}p_K}{ p_{23} p_{22}\cdots  p_{24-u}}
 < \frac{UB}{b_1}
\]
and the smallest $J$ verifying
\[
 \frac{ p_{24} p_{25}\cdots  p_{23+u}}{ p_{23} p_{22}\cdots  p_{25-u}p_J}
 < \frac{UB}{b_1}.
\]
After that we determined all integers 
$24\le k_1<\cdots < k_u \le K$ and
$23\ge j_1 > \cdots >j_u\ge J$ such that
\[
 \frac{ p_{k_1} p_{k_2}\cdots  p_{k_u}}{ p_{j_1} p_{j_2}\cdots  p_{j_u}}
 < \frac{UB}{b_1}.
\]
Each such tuple $(k_1, \ldots, k_u,j_1,\ldots, j_u)$ gives rise to
\[
\left\lfloor \frac{UB p_{j_1} p_{j_2}\cdots  p_{j_u}}{ b_1p_{k_1} p_{k_2}\cdots  p_{k_u}} \right\rfloor 
\]
candidates for the largest entry in a Diophantine couple $\{a,b\}$. 

Since the bound  $UD=10^{67.859}$ found in case (A) entails
$UB(a,UD)< 10^{33.9295} (16a^2+4a)^{-1/2}$, for $1\le a\le 7$ 
one has
\[
 \frac{UB(a,UD)}{b_1(a)} \le  \frac{UB(4,UD)}{b_1(4)} < 2.
\] 
Therefore, the multiplier $q$ mentioned above must be equal
to 1.

For each value of $b$ identified using the above method, 
we are able to show easily that there is no corresponding quadruple. 
This shows that $\omega(b) \leq 22$. In theory there is nothing 
stopping us from playing this trick again. However, when we search 
for $\omega(b) = 22$ we find that we could have over four thousand 
primes dividing $b$. This appears to be orders of magnitude harder 
than the $\omega(b) = 23$ case.

\subsection{Bounding $b$ in different ranges}\label{sec:6.2}
We have $ab +1 = r^2$. Note first that $d(r^{2} -1)$ is even (it is odd if and only if $r^{2} - 1 = s^{2}$ which implies that $(r+s)(r-s) = 1$ --- a contradiction). Since $d(r^{2}-1)$ counts the number of divisors of $r^{2} -1$ it follows that $\frac{1}{2}d(r^{2}-1)$ counts the number of \textit{pairs} of divisors $\{a, b\}$ with $a<b$. Now each $a$ corresponds with exactly one $b$ (and hence one pair corresponds with exactly one value of $a$): therefore $\frac{1}{2}d(r^{2}-1)$ is actually counting the divisors $a$. Furthermore, note that 
\begin{equation}\label{sum}
r^{2} -1 = ab > a^{2}.
\end{equation}
 Therefore $\frac{1}{2}d(r^{2}-1)$ is actually counting all those 
 $a$ with $a< \sqrt{r^{2} -1}$. Whence for a fixed $r$ we wish to count 
\[
\frac{1}{2} d_{\sqrt{r^{2} -1}}(r^{2} -1).
\]
If $r\leq R$ then summing over $r$ shows that the number of pairs $\{a, b\}$ is at most
\begin{equation}\label{legs}
\frac{1}{2} \sum_{r=3}^{R} d_{\sqrt{r^{2} -1}}(r^{2} -1) \leq \frac{1}{2} \sum_{r=3}^{R} d_{r}(r^{2} -1) <\frac{1}{2}\sum_{r=3}^{R} d_{R}(r^{2} -1).
\end{equation}

Now, we can make a slight improvement on (\ref{legs}). Since, for 
case (A) quadruples we have $b>4a$, we can improve on (\ref{sum}) 
to show that $r^{2} -1 = ab > 4a^{2}$.
Therefore, we amend (\ref{legs}) to show that the total number of 
pairs is at most
\begin{equation*}\label{legs2}
\frac{1}{2} \sum_{r=3}^{R} d_{R/2} (r^{2} -1).
\end{equation*}

One can go further than this. Let $N(\alpha, \beta)$ be the number of 
quintuples with $\alpha a < b \leq \beta a$, for some 
$\beta > \alpha \geq 4$. It then follows that for integers $m_{i}$ 
satisfying $4= m_{0}< m_{1} < \ldots < m_{k}$ the total number of 
quintuples is bounded above by
\begin{equation*}\label{fruit}
N(4, m_{1}) + N(m_{1}, m_{2}) + \cdots + N(m_{k-1}, m_{k}) + N(m_{k}, \infty),
\end{equation*}
where $N(m_{k}, \infty)$ means all those pairs $\{a, b\}$ such that 
$b> m_{k} a$. With the exception of $N(m_{k}, \infty)$, each number is 
of the form $N(m_{j}, m_{j+1})$.

Take $m_{j}a <b \leq m_{j+1} a$. Since 
$d> 4ab(4ab + a + b) > 16a^{2} b^{2} > 16 b^{4}/(m_{j+1})^{2}$ we have
\begin{equation}\label{omega}
b< \frac{d^{1/4} (m_{j+1})^{1/2}}{2}.
\end{equation}
We also have
\begin{equation}\label{sum3}
r^{2} -1 = ab > m_{j}a^{2} \Rightarrow a < R/\sqrt{m_{j}}.
\end{equation}
By taking $m_{j}$ large we ensure that the bound on $a$ in  
(\ref{sum3}) is small.
We now look at $\omega(b)$ for $b$ satisfying (\ref{omega}). 
We want to choose $m_{j+1}$ to be as large as possible such that 
we do not increase $\omega(b)$. For example, when $j=0$ we are 
considering $4a<b\leq m_{1}a$. We find, using $d \leq 7.228\times 10^{67}$, 
that $\omega(b)\leq 14$ provided that $m_{1} \leq 177$. Also, for 
$m_{2}$ we find that we can take $m_{2}\leq 499686$ and still ensure 
that $\omega(b)\leq 15$. We continue in this way, contenting ourselves 
with estimates on $m_{j}$ that are accurate to one decimal place. 
We find, using \textit{Mathematica} \cite{math}, that we may take
$$(m_{3}, m_{4}, m_{5}, m_{6}, m_{7}, m_{8}) = (1.7\cdot 10^{9}, 6.4\cdot 10^{12}, 2.9\cdot 10^{16}, 1.4\cdot 10^{20}, 7.8\cdot 10^{23},  4.8\cdot 10^{27}).$$
We know, from \S\ref{sec:6.1}, that there are at most 22 distinct prime factors of $b$.
Therefore we have that the number of quintuples is at most
\begin{equation*}
3 \cdot 2\left( 2^{15} \sum_{r=3}^{R} d_{R/2} (r^{2} -1) + 2^{16} \sum_{r=3}^{R} d_{R/\sqrt{177}}(r^{2} -1) + \cdots + 2^{23}\sum_{r=3}^{R} d_{R/\sqrt{4.8\cdot 10^{27}}}(r^{2} -1)\right).
\end{equation*}
We find that the above is no more than
\begin{equation}\label{goat}
1.177\cdot 10^{27}.
\end{equation}
 Using (\ref{case:b}), (\ref{case:c}), (\ref{case:d}) and (\ref{goat}) 
 we complete the proof of Theorem \ref{prince}.

Simion Stoilow Institute of Mathematics of the Romanian Academy, 
Research unit nr. 5, P.O. Box 1-764, RO-014700 Bucharest, Romania \\
Email: Mihai.Cipu@imar.ro \\[6pt]
Mathematical Sciences Institute, The Australian
National University, Canberra, Australia \\
Email: Timothy.Trudgian@anu.edu.au \\[6pt]

\end{document}